\documentclass[11pt]{amsart}

\usepackage{xypic}
\xyoption{all}

\newcommand{\myendbibitem}{\relax}
 \usepackage[linktocpage=true,hypertex]{hyperref}

\numberwithin{equation}{section}

\swapnumbers
\newtheorem{thm}[equation]{Theorem}
\newtheorem{prop}[equation]{Proposition}
\newtheorem{lem}[equation]{Lemma}
\newtheorem{cor}[equation]{Corollary}

\theoremstyle{definition}

\newtheorem{remark}[equation]{Remark}

\newcommand{\lrs}{\longrightarrow}
\newcommand{\Gr}{\operatorname{Gr}}
\newcommand{\Span}{\operatorname{Span}}
\newcommand{\mult}{\operatorname{mult}}

\newcommand{\bbG}{{\mathbb G}}
\newcommand{\bbC}{{\mathbb C}}

\newcommand{\bbZ}{{\mathbb Z}}

\newcommand{\Morph}{\operatorname{Morph}}
\newcommand{\Tmn}{T_{m,n}}

\newcommand{\Gmn}{G_{m,n}}

\newcommand{\inv}{^{-1}}
\newcommand{\Mnm}{(\Mn)^m}

\newcommand{\Z}{\operatorname{Z}}

\newcommand{\GL}{{\operatorname{GL}}}
\newcommand{\gl}{{\operatorname{\textit{gl}}}}
\newcommand{\SL}{{\operatorname{SL}}}
\newcommand{\liesl}{{\operatorname{\textit{sl}}}}
\newcommand{\PGL}{{\operatorname{PGL}}}
\newcommand{\PGLn}{{\operatorname{PGL}_n}}

\newcommand{\Stab}{\operatorname{Stab}}

\newcommand{\Char}{\operatorname{char}} 
\newcommand{\diag}{\operatorname{diag}}

\newcommand{\RMaps}{\operatorname{RMaps}}
\newcommand{\lra}{\longrightarrow}
\newcommand{\Mat}{{\operatorname{M}}}
\newcommand{\M}{\operatorname{M}}
\newcommand{\Mn}{\Mat_n}

\newcommand{\trdeg}{\operatorname{trdeg}}
\newcommand{\Sym}{{\operatorname{S}}}
\newcommand{\tr}{\operatorname{tr}}
\newcommand{\pr}{\operatorname{pr}}

\newcommand{\SLm}{{\SL_m}}
\newcommand{\GLm}{{\GL_m}}

\newcommand{\UD}{\operatorname{{\it UD}}}
\newcommand{\UDmn}{\UD(m,n)}
\newcommand{\Lie}{\operatorname{\it Lie}}
\DeclareMathOperator{\catq}{/\!/}

\begin{document}

\date{September 2, 2005}

\title[Group actions and invariants]
   {Group actions and invariants in algebras of generic matrices}

\author{Z. Reichstein}
\address{Department of Mathematics, University of British Columbia,
       Vancouver, BC V6T 1Z2, Canada}
\email{reichste@math.ubc.ca}
\urladdr{www.math.ubc.ca/$\stackrel{\sim}{\phantom{.}}$reichst}
\thanks{Z. Reichstein was supported in part by an NSERC research
  grant}

\author{N. Vonessen}
\address{Department of Mathematical sciences, University of Montana,
  Missoula, MT 59812-0864, USA}
\email{Nikolaus.Vonessen@umontana.edu}
\urladdr{www.math.umt.edu/vonessen}
\thanks{N.\ Vonessen gratefully acknowledges the support of the
  University of Montana and the hospitality of the University of
  British Columbia during his sabbatical in 2002/2003, when part of
  this research was done.}




\begin{abstract}
  We show that the fixed
  elements for the natural $\GL_m$-action on the universal division
  algebra $\UD(m, n)$ of $m$ generic $n \times n$-matrices form a
  division subalgebra of degree $n$, assuming $n\geq 3$ and $2 \le m
  \le n^2 -2$. This allows us to describe the asymptotic behavior 
  of the dimension of the space of $\SLm$-invariant homogeneous central 
  polynomials $p(X_1, \dots, X_m)$ for $n \times n$-matrices.
  Here the base field is assumed to be of characteristic zero. 
\end{abstract}

\maketitle

\tableofcontents

\section{Introduction}
\label{sect.intro}

Let $k$ be a field of characteristic zero, $m$ and $n$ be integers
$\ge 2$, and $\Gmn = k \{ X_1, \dots, X_m \}$ be the $k$-algebra of
$m$ generic $n \times n$-matrices. That is, $\Gmn$ is the
$k$-subalgebra of $\M_n(k[x_{ij}^{(h)}])$ generated by
\[ X_1 = (x_{ij}^{(1)}) , \dots, X_m = (x_{ij}^{(m)}) \, , \]
where the $x_{ij}^{(h)}$ are $m n^2$ independent commuting variables.
By a theorem of Amitsur, $\Gmn$ 
is a domain of PI-degree $n$.
There is a natural action of the general linear group $\GLm$ on $\Gmn$
given by
\begin{equation}\label{GLm-action}
 g = (g_{ij}) \colon X_j \longmapsto \sum_{i = 1}^m g_{ij} X_i \, .
\end{equation}
In this paper we will prove the following theorem.
 
\begin{thm}\label{thm.intro1}
For $2\leq m\leq n^2-2$, the domain $(\Gmn)^\SLm$ has PI-degree~$n$.
\end{thm}

The trace ring $\Tmn$ of $\Gmn$ is defined as the subring 
of $\M_n(k[x_{ij}^{(h)}])$ generated by elements of
the form $Y$ and $\tr(Y)$, as $Y$ ranges over $\Gmn$.
The action~\eqref{GLm-action} on $\Gmn$ naturally extends to
$\Tmn$.  Note that the algebras $\Gmn$ and $\Tmn$, and their centers
$\Z(\Gmn)$ and $\Z(\Tmn)$ have a natural $\bbZ$-grading
inherited from $\M_n(k[x_{ij}^{(h)}])$ 
(each variable $x_{ij}^{(h)}$ has degree $1$) and that this grading
is preserved by the action~\eqref{GLm-action}. 
As a consequence of Theorem~\ref{thm.intro1} we obtain the following
result. 

\begin{thm} \label{thm.intro2}
  Let $2\leq m\leq n^2-2$, and let $R$ be one of the rings
  $\Gmn$, $\Tmn$, $\Z(\Gmn)$, or $\Z(\Tmn)$. Denote the degree $d$
  homogeneous component of $R$ by $R[d]$. Then
  \[\limsup_{d\to\infty}\,\frac{\dim_k \, R^{\SLm}[d]}%
             {d^{(m-1)n^2 - m^2 + 1}}\]
  is a finite nonzero number.
\end{thm}  

One can think of the center of $\Gmn$ as consisting of the
$m$-variable central polynomials for $n\times n$-matrices (over
commutative $k$-algebras).  Theorem~\ref{thm.intro2} thus describes, for
$R=\Z(\Gmn)$, the asymptotic behavior of the dimension of the space of
$\SLm$-invariant homogeneous central polynomials $p(X_1, \dots, X_m)$
for $n \times n$-matrices.

The $\GLm$-representations on $\Gmn$, $\Z(\Gmn)$, $\Tmn$ and $\Z(\Tmn)$ 
have been extensively studied; see, 
e.g.,~\cite{berele1, berele2, drensky, formanek3, regev}. 
Once again, let $R$ be one of these rings.
Recall that the irreducible polynomial representations of
$\GLm$ are indexed by partitions $\lambda = (\lambda_1, \dots, \lambda_s)$
with~$s \le m$ parts; cf., e.g.,~\cite[Section 2]{formanek3}. 
Denote the multiplicity of the irreducible $\GLm$-representation 
corresponding to~$\lambda$ in $R$ by 
$\mult_{\lambda}(R)$. If $(r^m)$ is the partition $(r, \dots, r)$ 
($m$ times) then it is easy to show that
\[ \dim \, R^{\SLm}[d] = \begin{cases} 
\,\mult_{(r^m)}(R) &\text{if $d = rm$,} \\
\,0 &\text{if $d$ is not a multiple of $m$;} \end{cases}
\] 
see Remark~\ref{rem9.2}. The conclusion of 
Theorem~\ref{thm.intro2} can thus be rephrased by saying that
  \[\limsup_{r\to\infty}\,\frac{\mult_{(r^m)}(R)} {r^{(m-1)n^2 - m^2 + 1}}\]
is a finite nonzero number.
We also note that by the Berele-Drensky-Formanek 
correspondence, $\mult_{(r^m)}(R)$ equals the multiplicity 
of the $S_m$-character $\chi^{(d^m)}$ in the cocharacter 
sequence of $R$; see~\cite[Section 4]{formanek3}. 

\smallskip
The division algebra of quotients of $\Gmn$ (or equivalently, 
of $\Tmn$) is called the {\em universal division algebra} 
of $m$ generic $n \times n$-matrices; we shall denote it by $\UDmn$.
Note that the $\GLm$-action~\eqref{GLm-action} on $\Gmn$
naturally extends to $\UDmn$. We shall deduce Theorem~\ref{thm.intro1}
from the following related result.

\begin{thm}\label{thm:GLm}
  If $2\leq m\leq n^2-2$ and $n\geq3$, then $\UDmn^\GLm$ is a division
  algebra of degree n. 
\end{thm}

For all other values of $m, n \ge 2$, $\UDmn^\GLm$ is a field; see
Propositions~\ref{prop.m=n^2} and~\ref{prop.m=n^2-1}. A brief 
summary of the properties of $\UDmn^{\GLm}$ and $\UDmn^{\SLm}$ 
is given in the two tables below. 

\bigskip

\newsavebox{\scratchbox}
\savebox{\scratchbox}{\footnotesize $m\leq n^2-2$, $n\geq3$}

\begin{table}[ht]
\footnotesize   
\caption{\strut Properties of $\UDmn^\GLm$}
\newcommand{\mystrut}{\rule{0cm}{4mm}}    
\begin{tabular}{l|c|c|c}
Case \mystrut&   PI-Degree & Transcendence Degree/$k$ & Central in $\UDmn$? \\ \hline
\usebox{\scratchbox}& $n$\mystrut& $(m-1)n^2-m^2+1$ & No \\ 
$m=n^2-1$    & 1 & $n-1$ & No \\
$m=n=2$  & 1 & 1 & No     \\
$m\geq n^2$ & 1  & $0$  & Yes \\ 
\end{tabular}                             
\end{table}

\begin{table}[ht]
\footnotesize   
\caption{\strut Properties of $\UDmn^\SLm$}
\newcommand{\mystrut}{\rule{0cm}{4mm}}    
\begin{tabular}{l|c|c|c}
Case \mystrut&   PI-Degree & Transcendence Degree/$k$ & Central in $\UDmn$? \\ \hline
\makebox[0pt][l]{$m\leq n^2-2$}\hphantom{\usebox{\scratchbox}} 
    & $n$\mystrut& $(m-1)n^2-m^2+ 2$ & No \\ 
$m=n^2-1$    & 1 & $n$ & No \\
$m=n^2$  & 1 & 1 & Yes     \\
$m\geq n^2+1$ & 1  & $0$  & Yes \\ 
\end{tabular}                             
\end{table}  

\bigskip
\enlargethispage{-5mm}

The assertions of the tables in the cases where $m \le n^2 -2$ and
$n\geq3$ are based on
Theorems~\ref{thm:GLm} and~\ref{thm1.grass}, the case where $m = n = 2$ 
is considered in~\cite[Section 14]{rv3}, and the cases 
where $m \ge n^2 -1$ are treated in detail in Section~\ref{sect.m>=n^2-1}.

\smallskip
It appears likely that Theorems~\ref{thm.intro1} --~\ref{thm:GLm} 
remain valid in prime characteristic (perhaps, not 
dividing $n$); we have not attempted to extend them
in this direction. Our arguments rely on the work 
of Richardson~\cite{richardson} and on our own prior papers
\cite{rv2,rv3}, all of which make the characteristic zero assumption.%
\footnote{We remark that Richardson~\cite{richardson} worked
over $k = \bbC$, and his proofs are based on analytic 
techniques. The results we need (in particular,
\cite[Theorem 9.3.1]{richardson}), remain valid over any 
algebraically closed field of characteristic 
zero by the Lefschetz principle. 
Extending~\cite[Theorem 9.3.1]{richardson} to 
prime characteristic is an open problem of independent interest.}
 
\subsection*{Conventions and Terminology.}
All central simple algebras in this paper are assumed to be
finite-dimensional over their centers.  All algebraic varieties,
algebraic groups, group actions, morphisms, rational maps, etc., are
assumed to be defined over the base field $k$ (which we always assume
to be of characteristic zero).  By a point of an algebraic variety $X$
we shall always mean a $k$-point.
Throughout, $G$ will denote a linear algebraic group.  We shall refer to
an algebraic variety $X$ endowed with a regular $G$-action as a
$G$-variety.  We will say that a $G$-variety $X$ (or the $G$-action 
on $X$) is generically free if $\Stab_G(x) = \{ 1 \}$ for $x \in X$ 
in general position.
Finally, unless otherwise specified, $m$ and $n$ are integers $\geq 2$.

\section{Preliminaries}
\label{sect.geometric-actions}

\subsection*{Concomitants}

Let $\Gamma$ be an algebraic group and $V$ and $W$ be $\Gamma$-varieties.
Then we shall denote the set of $\Gamma$-equivariant morphisms $V \lra W$
(also known as {\em concomitants}) by $\Morph_{\Gamma}(V, W)$ 
and the set of $\Gamma$-equivariant rational maps $V \dasharrow W$ 
(also known as {\em rational concomitants}) by $\RMaps_{\Gamma}(V, W)$.

In the case where $W$ is a finite-dimensional linear representation
of $\Gamma$, we also define a {\em relative concomitant} as a 
morphism $f \colon V \lra W$ satisfying the following condition
(which is slightly weaker than $\Gamma$-equivariance):
there is a character $\chi \colon \Gamma \lra k^*$ such that
\[ f(g \cdot v) = \chi(g) \ \bigl(g \cdot f(v)\bigr) \]  
for all $v \in V$ and $g \in \Gamma$. 
For a rational map $f \colon V \dasharrow W$ 
the notion of a {\em relative rational concomitant} is defined 
in a similar manner. 
If $W = k$, with trivial $\Gamma$-action, then the term ``invariant" 
is used in place of ``concomitant". For future reference we record 
the following: 

\begin{lem} \label{lem.dc} Suppose $V$ and $W$ are finite-dimensional
linear representations of $\Gamma$. Every rational concomitant 
$f \colon V \dasharrow W$ can be written as $\frac{a}{b}$,
where $a$ is a relative concomitant and $b$ is a relative 
invariant.
\end{lem}

\begin{proof} See the proof of~\cite[Chapter~1, Proposition~1]{dc}.
  Note that the characters associated to $a$ and $b$ are necessarily
  equal. 
\end{proof}

If $W$ is a $k$-algebra and $\Gamma$ acts on $W$ by $k$-algebra 
automorphisms, then the algebra structure of $W$ 
induces algebra structures on $\Morph_{\Gamma}(V, W)$ and
$\RMaps_{\Gamma}(V, W)$ in a natural way.  Namely, given $a, b \colon V \lra W$
(or $a, b \colon V \dasharrow W$), one defines $a + b$ and $ab$ 
by $(a+b)(v) = a(v) + b(v)$ and $ab(v) = a(v) b(v)$ for $v \in V$.

\begin{thm} \label{thm.procesi}
{\em (Procesi~\cite[Theorem 2.1]{procesi2}; cf.\ also
\cite[Theorem 10]{formanek3}, or \cite[Theorem 14.16]{saltman})}
Let $\Mnm$ be
the space of $m$-tuples of $n\times n$-matrices;
the group $\PGLn$ acts on it by simultaneous conjugation.  Then
\begin{itemize}
\item[(a)] $\Morph_{\PGLn}(\Mnm, \Mn) \simeq \Tmn$
\item[(b)] $\RMaps_{\PGLn}(\Mnm, \Mn) \simeq \UDmn$
\end{itemize}
Moreover, the two isomorphisms identify the $i$-th projection
$\Mnm \lra \Mn$ with the $i$-th generic matrix $X_i$. 
\qed
\end{thm}

Here $\Tmn$ and $\UDmn$ are, respectively, the trace ring and the
universal division algebra of $m$ generic $n \times n$-matrices,
defined in the introduction.  Note that part~(b) of
Theorem~\ref{thm.procesi} follows from part~(a) by Lemma~\ref{lem.dc},
since the simple group $\PGLn$ does not have nontrivial characters
(so that relative concomitants and invariants are actually
concomitants and invariants, respectively). 

We also remark that the construction of $\Tmn$ remains well-defined 
if $m=1$.  Theorem~\ref{thm.procesi} also holds in this case, provided 
that one defines $\UD(1,n)$ to be the field of quotients of $T_{1,n}$, 
rather than $G_{1,n}$. (For $m \ge 2$, $\Tmn$ and $\Gmn$ have 
the same division algebra of quotients, but this is not 
the case for $m = 1$.).

\subsection*{Geometric actions}
For the rest of this section we will assume that $k$ is algebraically 
closed.  If $X$ is a $\PGLn$-variety, then, as we mentioned above,
$\RMaps_{\PGLn}(X, \Mn)$ has an algebra structure
naturally induced from $\Mn$.  If the $\PGLn$-action 
on $X$ is generically free then $\RMaps_{\PGLn}(X,
\Mn)$ is a central simple algebra of degree $n$, with center
$k(X)^{\PGLn}$; cf.\ \cite[Lemmas 8.5 and 9.1]{reichstein}.  

Suppose that $X$ is a $G \times \PGLn$-variety, and that the
$\PGLn$-action on $X$ is generically free. Then the $G$-action on $X$
naturally induces a $G$-action on $\RMaps_{\PGLn}(X, \Mn)$. 
Following~\cite{rv3} we define the action of an algebraic group $G$ on a
central simple algebra $A$ to be {\em geometric} if $A$ is
$G$-equivariantly isomorphic to $\RMaps_{\PGLn}(X, \Mn)$ for some $G \times
\PGLn$-variety $X$ as above. The $G \times \PGLn$-variety $X$ is then
called the {\em associated variety} for the $G$-action on $A$; this
associated variety is unique (as a $G \times \PGLn$-variety), up to
birational isomorphism; cf. \cite[Corollary 3.2]{rv3}.

Note that we defined geometric actions only if $k$ is algebraically
closed.  Also note that if an algebraic group acts geometrically on a
central simple algebra $A$, then the center of $A$ is necessarily a
finitely generated field extension of $k$.

Of particular interest to us will be the case where 
$X =\Mnm$ is the space of $m$-tuples of $n\times n$-matrices.  Here
$\PGLn$ acts on $\Mnm$ by simultaneous conjugation (since $m \ge 2$,
this action is generically free) and $G = \GLm$ acts 
on $(A_1, \dots, A_m)\in\Mnm$ by sending $(A_1, \dots, A_m)$
to $(B_1, \dots, B_m)$ where $B_j = \sum_{i = 1}^m
c_{ij} A_i$ and $g^{-1} = (c_{ij})$.  The actions of $\GLm$
and $\PGLn$ commute, and the $\GLm$-action on $\Mnm$ induces the
$\GLm$-action~\eqref{GLm-action} on $\UDmn$.  So $\Mnm$ is the
associated variety for the $\GLm$-action on $\UDmn$; see
Theorem~\ref{thm.procesi} (cf.\ also \cite[Example 3.4]{rv3}).

We conclude this section with a simple result which we will use
repeatedly. 

\begin{lem} \label{lem.stab}
  Assume $k$ is algebraically closed.  Let $X$ be a
  $G\times\PGLn$-variety which is $\PGLn$-generically free.  Denote by
  $\pi\colon X\dasharrow X/G$ the rational quotient map for the
  $G$-action.  Then for $x\in X$ in general position, the projection
  $\pr_2 \colon G \times \PGLn \lra \PGLn$ onto the second factor induces an
  isomorphism from $\Stab_{G \times \PGLn}(x)/\Stab_G(x)$ onto
  $\Stab_{\PGLn}(\pi(x))$.  \qed
\end{lem}

\begin{proof} Recall that by a theorem of Rosenlicht, $\pi^{-1}(\overline{x})$
is a single $G$-orbit for $\overline{x} \in X/G$ in general position; 
see~\cite[Theorem 2]{rosenlicht} or \cite[Section 2.3]{pv}. 
Consequently, for $x \in X$ in general position
the projection $\pr_2$ restricts to a surjective morphism 
$\Stab_{G \times \PGLn}(x) \lra \Stab_{\PGLn}(\pi(x))$
of algebraic groups.  The kernel of this morphism is clearly  
$\Stab_G(x)$, and the lemma follows.
\end{proof}

\section{Geometric actions on division algebras} 
\label{sect.div}

Throughout this section we will assume that $k$ is algebraically closed.
The main result of this section is the following theorem; after
its proof, we will deduce several corollaries.

\enlargethispage{4mm}

\begin{thm}\label{thm.div}
  Assume $k$ is algebraically closed.
  Let $G$ be a linear algebraic group acting geometrically on a
  division algebra $D$ of degree~$n$.  Let $X$ be the associated
  $G\times\PGLn$-variety.  Then for $x\in X$ in general position,
  \[S_x:=\Stab_{G \times \PGL_n}(x)/ \Stab_G(x)\]
  is reductive.
\end{thm}

\begin{proof}
Let $X$ be the associated $G \times \PGLn$-variety for the $G$-action
on $D$. Recall that the $\PGLn$-action on $X$ is generically free.
We want to show that the group
$S_x = \Stab_{G \times \PGLn}(x)/\Stab_G(x)$ is 
reductive for $x \in X$ in general position.  Assume the contrary.
Denoting the unipotent radical by $R_u$, this means that 
$R_u(\Stab_{G \times \PGLn}(x))$ is not contained in $G$.
Since unipotent groups are connected, this is equivalent to 
\begin{equation} \label{e.contr}
\Lie\bigl(R_u(\Stab_{G \times \PGLn}(x))\bigr) \not \subseteq \Lie(G) 
\end{equation}
for $x \in X$ in general position. Here and in the sequel 
$\Lie$ stands for the Lie algebra.  To simplify notation, set
$H=G\times \PGLn$, and for $x\in X$, set $H_x=R_u(\Stab_H(x))$.
Now define $U_X\subseteq X \times \Lie(H)$ by
\[ U_X = \bigl\{ (x, a) \, \bigl| \, x\in X \text{ and } a \in 
\Lie(H_x)\bigr\}\, . \]
We first show that $U_X$ is a vector bundle over a dense open subset
$X_0 \subset X$.  By~\cite[6.2.1, 9.2.1, and 6.5.3]{richardson}, there
is an $H$-stable dense open subset $X_0$ of $X$ such that $\{H_x\mid
x\in X_0\}$ is an algebraic family of algebraic subgroups of $H$.
Moreover, $\dim(H_x)$ is constant for $x\in X_0$, say equal to~$d$.
Replacing $X$ by $X_0$, we may assume that $\{H_x\mid x\in X\}$ is an
algebraic family of algebraic subgroups of $H$.
By~\cite[6.2.2]{richardson}, $x\mapsto\Lie(H_x)$ defines a morphism of
algebraic varieties from $X$ to the Grassmannian
of $d$-dimensional subspaces of $\Lie(H)$.  Since the universal bundle
over this Grassmannian is a vector bundle (see, e.g.,
\cite[3.3.1]{weyman}), its pull-back $U_X$ is a vector bundle over $X$.

Note also that $U_X$ is, by definition, an $H$-invariant subbundle of
the trivial bundle $X \times \Lie(H) \lra X$; here $H$ acts on its Lie
algebra by the adjoint action.  Since the $\PGLn$-action on $X$ is
generically free, the no-name lemma tells us that there is a
$\PGL_n$-equivariant birational isomorphism $U_X \dasharrow X \times
k^d$ such that the following diagram commutes
\[ \xymatrix{
U_X \ar@{->}[d]  \ar@{-->}[r]^-{\simeq\,\,} & X \times k^d \ar@{->}[ld] \cr
X  & } \]
(For a proof and a brief discussion of the history of the no-name lemma,
see~\cite[Section 4.3]{cgr}.)  In other words, the vector
bundle $U_X \lra X$ has $d$ $\PGLn$-equivariant rational sections
$\beta_1, \dots, \beta_d \colon X \dasharrow U_X$ such that
$\beta_1(x)$, \ldots, $\beta_d(x)$ are linearly independent for $x \in X$
in general position.  We identify here $\beta_i(x)$ with $a$ if
$\beta_i(x)=(x,a)\in \{x\}\times\Lie(H_x)$.  In view
of~\eqref{e.contr}, some $k$-linear combination $\beta = c_1 \beta_1 +
\dots + c_d \beta_d$ has the property that $\beta(x) \not \in \Lie(G)$
for $x \in X$ in general position.

Now recall that the natural projection $\SL_n\lra\PGLn$ induces
a Lie algebra isomorphism $\liesl_n\lra\Lie(\PGLn)$, allowing us to
identify the two Lie algebras.  Hence
\[ U_X \subseteq  X\times \Lie(G) \times \liesl_n \, . \]
Let $f = pr \circ \beta \colon X \dasharrow \liesl_n$, where $pr \colon
U_X \lrs \liesl_n$ denotes the natural projection.  Note that $\liesl_n
\hookrightarrow \gl_n=\M_n$, so that $f$ may be viewed as a
$\PGLn$-equivariant rational map $X \dasharrow \Mn$, i.e., as an
element of $D$.  The condition that $\beta(x) \not \in \Lie(G)$
ensures that $f \neq 0$.  On the other hand, we will show below that
for $x \in X$ in general position, $f(x)$ is a nilpotent $n \times
n$-matrix, so that $f^n = 0$. This means that $D$ contains a non-zero
nilpotent element~$f$, contradicting our assumption that $D$ is a
division algebra.

It remains to be shown that for $x\in X$ in general position, $f(x)$
is a nilpotent matrix.  The natural projection $G\times\PGLn\lra\PGLn$
maps the unipotent group $H_x$ to a unipotent subgroup~$U$ of $\PGLn$.
Denote by $K$ the preimage of $U$ in $\SL_n$.  It is a solvable group,
so its subset $K_u$ of unipotent elements is a closed subgroup.  The
surjection $K_u\lra U$ is finite-to-one, so their Lie algebras are
isomorphic.  In particular, $f(x)$ belongs to
$\Lie(K_u)\subset\liesl_n\subset\gl_n=\Mn$.  Finally, since $K_u$ is
a unipotent subgroup of $\GL_n$, its Lie algebra in $\Mn$ consists of
nilpotent matrices, see, e.g., \cite[I.4.8]{borel}.  This completes
the proof of Theorem~\ref{thm.div}.
\end{proof}

We now proceed with the corollaries. Recall that a subgroup 
$S \subset \Gamma$ is said to be a {\em stabilizer in general 
position} for a $\Gamma$-variety $X$ if there exists a dense
$\Gamma$-invariant subset $U \subset X$ such that $\Stab(x)$
is conjugate to $S$ for any $x \in U$. For a detailed discussion 
of this notion, see~\cite[Section 7]{pv}. 

\begin{cor}\label{thm.div.cor1}
  Assume $k$ is algebraically closed.
  Let $G$ be a linear algebraic group acting geometrically on a
  division algebra $D$ of degree~$n$.  Let $X$ be the associated
  $G\times\PGLn$-variety.
  \begin{itemize}
  \item[(a)]The induced $\PGLn$-action on the rational quotient $X/G$
    has a stabilizer~$S$ in general position.  Moreover, $S$ is
    reductive, and $S\simeq S_x=\Stab_{G \times \PGL_n}(x)/
    \Stab_G(x)$ for $x\in X$ in general position.
  \item[(b)]If the $G$-action on $X$ is generically free, then
  \begin{align*}
    \trdeg_k(\Z(D^G)) &= \trdeg_k(\Z(D)^G)\\
                     &= \dim(X)-\dim(G)+\dim(S) - n^2 + 1 \,.
  \end{align*}
  \end{itemize}
\end{cor}

\begin{proof}
  (a) It follows from Theorem~\ref{thm.div} and Lemma~\ref{lem.stab}
  that points in $X/G$ in general position have a reductive
  stabilizer.  A theorem of Richardson (see~\cite[Theorem
  9.3.1]{richardson} or \cite[Theorem 7.1]{pv}) now implies that the
  $\PGLn$-action on $X/G$ has a stabilizer~$S$ in general position.
  By Lemma~\ref{lem.stab}, $S\simeq S_x=\Stab_{G \times \PGL_n}(x)/
  \Stab_G(x)$ for $x\in X$ in general position.

  \smallskip 
  
  (b) The first equality follows from the fact that $\Z(D^G)$ is an
  algebraic extension of $\Z(D)^G$. Indeed, the minimal polynomial of
  any element of $D^G$ over $\Z(D)$ is unique and must therefore have
  coefficients in $\Z(D)^G$.
  
  To prove the second equality, note that, $\Z(D) = k(X/\PGLn) =
  k(X)^{\PGLn}$ and thus
  \[ \Z(D)^G \simeq (k(X)^{\PGLn})^G = k(X)^{G \times \PGLn} =  
  k(X/(G \times \PGLn)) \, . \]
  Since we are assuming that $G$ acts generically freely on $X$,
  part~(a) implies that $S\simeq \Stab_{G\times\PGLn}(x)$ for $x\in X$
  in general position.  Hence the dimension of the general fiber of
  the rational quotient map $X \dasharrow X/(G \times \PGLn)$ is equal
  to the dimension of $(G \times \PGLn)/S$.  The fiber dimension
  theorem now tells us that the transcendence degree of $\Z(D)^G$ is
  \[ \dim \, X/(G \times \PGLn) = \dim(X) - \dim(G) - \dim(\PGLn)
      + \dim(S)\,. \qedhere \]
\end{proof}

\begin{cor}
  Assume $k$ is algebraically closed.  Let $G$ be a unipotent linear
  algebraic group acting geometrically on a division algebra $D$ of
  degree~$n$. Then $D^G$ is a division algebra of degree $n$.
\end{cor}

This was proved for algebraic actions in \cite[Proposition~12.1]{rv3}.

\begin{proof}
  By \cite[Lemma~7.1]{rv3}, for $x\in X$ in general position,
  $\Stab_{G \times \PGLn}(x)$ is isomorphic to a subgroup of $G$, so
  is unipotent. On the other hand, by Theorem~\ref{thm.div}, the
  projection $S_x$ of this group to $\PGLn$ is reductive.  Thus $S_x$
  is both unipotent and reductive, which is only possible if $S_x = \{
  1 \}$.  In other words,
  \[ \Stab_{G \times \PGLn}(x)  \subset G \times \{ 1 \} \, . \]
  The desired conclusion now follows from \cite[Theorem~1.4]{rv3}.
\end{proof}

\section{Dimension counting in the Grassmannian}
\label{sect.dim-counting}

In preparation for the proof of Theorem~\ref{thm:GLm}
in the next section, we will now establish the following: 

\begin{prop} \label{prop.grass1}
  Assume $k$ is algebraically closed.
  Let $V$ be an $N$-di\-men\-sion\-al $k$-vector space and let $V = V_1
  \oplus \dots \oplus V_r$, where $\dim(V_i) = N_i \ge 1$.  Let $Z$ be
  the subset of the Grassmannian $\Gr(m, N)$ consisting of
  $m$-dimensional subspaces $W$ of $V$ of the form
  \[ W = W_1 \oplus \dots \oplus W_r \, , \]
  where $W_i \subseteq V_i$. {\textup(}Here we allow $W_i = (0)$ for
  some $i$.\textup{)} Then $Z$ is a closed subvariety of $\Gr(m,
  N)$. If $2\leq m\leq N - 2$, then each irreducible component of $Z$
  has codimension $\ge N - \max_{i = 1, \dots, r}(N_i)$ in $\Gr(m,
  N)$.  Moreover, equality holds \textup{(}for some irreducible
  component of $Z$\textup{)} only if \textup{(i)} $r = 1$ or
  \textup{(ii)} $r = 2$, $m = 2$ and $N = 4$.
\end{prop}

\begin{proof}
Let $m_1, \dots m_r$ be non-negative integers such that $m_1 + \dots
+ m_r = m$ and such that $m_i\leq N_i$ for all $i$.  Let $Z_{m_1,
    \dots , m_r}$ be the image of the map
  \[ \phi_{m_1, \dots, m_r} \colon \Gr(m_1, N_1) \times \dots \times
     \Gr(m_r, N_r ) \lra \Gr(m,N)
  \]
  given by $(W_1, \dots , W_r) \mapsto W_1 \oplus \dots \oplus W_r$.
  (Here $\Gr(m_i, N_i)$ is the Grassmannian of $m_i$-dimensional
  vector subspaces of $V_i$.)  Since the domain of the map $\phi_{m_1,
    \dots, m_r}$ is projective, its image is closed in $\Gr(m, N)$.
  Thus each $Z_{m_1, \dots, m_r}$ is a closed
  irreducible subvariety of $\Gr(m, N)$ birationally isomorphic to
  \[ \Gr(m_1, N_1) \times \dots  \times \Gr(m_r, N_r) \]
  and $Z$ is the union of the $Z_{m_1, \dots, m_r}$.  It 
  remains to show that
  \begin{equation}\label{eqn:sec4:1}
      \dim \, \Gr(m, N) - \sum_{i=1}^r\dim \, \Gr(m_i, N_i)
      \ge N - \max_{i = 1, \dots, r} (N_i)  \, ,
  \end{equation}
and that equality is only possible if $r = 1$ or $r = 2$,  $N_1 = N_2 = 2$ and
$m_1 = m_2 = 1$ (and thus $N = N_1 + N_2 = 4$ and $m = m_1 + m_2 = 2$).
Recall that $\dim\Gr(m,N)=(N-m)m$.
Letting $l_i = N_i - m_i$ and $l = N-m = l_1 + \dots + l_r$,
we can rewrite~\eqref{eqn:sec4:1} as
\[ lm - \sum_{i = 1}^r l_i m_i \ge l + m -  
    \max_{i = 1, \dots, r} \, (l_i + m_i) \]
or, equivalently,
\[ (l - 1) (m-1) - 1 \ge \sum_{i = 1}^r l_i m_i - 
 \max_{i = 1, \dots, r} (l_i + m_i) \, . \]
Proposition~\ref{prop.grass1} is thus a consequence
of the following elementary lemma.
\end{proof}

\begin{lem} \label{lem.elem}
Let $(l_1, m_1), \dots, (l_r, m_r)$ be $r$ pairs of non-negative integers
and let $l = \sum_{i = 1}^r l_i$ and $m = \sum_{i = 1}^r m_i$. Assume that
$l_i + m_i \ge 1$ for every $i = 1, \dots, r$ and $l, m \ge 2$.
Then
\begin{equation} \label{e.elem}
(l - 1) (m-1) - 1 \ge \sum_{i = 1}^r l_i m_i - 
 \max_{i = 1, \dots, r} (l_i + m_i) \, .
\end{equation}
Moreover, equality holds if and only if either \textup{(i)} $r = 1$ or
\textup{(ii)} $r = 2$ and $(l_1, m_1) = (l_2, m_2) = (1, 1)$.
\end{lem}

\begin{proof}We consider two cases.

\smallskip
{\em Case 1:} Suppose that for every $i = 1, \dots, r$, 
either $l_i = 0$ or $m_i = 0$.  Since $l, m \ge 2$, 
we have
$ (l - 1) (m-1) - 1 \ge 0$. On the other hand,
$\sum_{i = 1}^r l_i m_i - 
 \max_{i = 1, \dots, r} (l_i + m_i) = 
- \max_{i = 1, \dots, r} (l_i + m_i)  < 0$. 
Hence, in this case~\eqref{e.elem} holds and is a strict 
inequality.

\smallskip
{\em Case 2:} Now suppose that $l_i , m_i \ge 1$ for some $i \ge 1, \dots, r$.
After renumbering the pairs $(l_1, m_1), \dots, (l_r, m_r)$, we may assume
$i = 1$. Now set
\[ l_j' = \begin{cases} \text{$l_1 - 1$, if $j = 1$,} \\
 \text{$l_j$, otherwise;} \end{cases} \quad \text{and} \quad 
 m_j' = \begin{cases} 
\text{$m_1 - 1$, if $j = 1$} \\
 \text{$m_j$, otherwise.} 
\end{cases} \] 
Note that $l_j', m_j' \ge 0$ for every $j = 1, \dots, r$. Thus
\begin{align*}
(l - 1) (m-1) - 1 &= (\sum_{i = 1}^r l_i' ) (\sum_{j = 1}^r m_j') - 1= 
\sum_{i = 1}^r l_i' m_i' + \sum_{i \ne j} l_i' m_j' - 1\\
& \ge\sum_{i = 1}^r l_i' m_i' - 1 = 
\sum_{i = 1}^r l_i m_i  - (l_1 + m_1)\\
& \ge  
\sum_{i = 1}^r l_i m_i  - \max_{i = 1, \dots, r} (l_i + m_i) \, . 
\end{align*}
This completes the proof of the inequality~\eqref{e.elem}.

\smallskip
It is easy to see that equality holds in cases (i) and (ii). 
It remains to show that the inequality~\eqref{e.elem}
is strict for all other choices of $(l_1, m_1), \dots, (l_r, m_r)$.
Indeed, a closer look at the above argument shows that equality 
in~\eqref{e.elem} can hold if and only if we are in Case 2 and

\smallskip
(a) $l_i' m_j' = 0$ whenever $i \ne j$ and 

\smallskip
(b) $l_1 + m_1 = \max_{i = 1, \dots, r} (l_i + m_i)$.

\smallskip
Assume that conditions (a) and (b) are satisfied.  Since 
$\sum_{i = 1}^r l_i' = l- 1 \ge 1$, we cannot have $l_i' = 0$ 
for all $i = 1, \dots, r$. In other words, $l_{i_0}' \ge 1$ for
some $i_0 \in \{ 1, \dots, r \}$. Then condition (a)
says that $m_j' = 0$ for every $j \ne i_0$. On the other hand,
$m_{i_0}' = \sum_{j = 1}^r m_j' = m- 1 \ge 1$, and applying
condition (a) once again, we conclude that
that $l_i' = 0$ for every $i \ne i_0$. To sum up, there exists
an $i_0 \in \{ 1, \dots, r \}$ such that $l_{i_0}' = l - 1$, 
$m_{i_0} = m-1$ and $l_i' = m_i' = 0$ for every $i \ne i_0$.

In particular, for every $i \ne 1, i_0$, we have
$l_i = l_i' = 0$ and $m_i = m_i' = 0$, contradicting our assumption
that $l_i + m_i \ge 1$. Thus $i \in \{  1, i_0 \}$ 
for every $i = 1, \dots, r$.
In other words, either $i_0 = 1$ and $r = 1$ 
(in which case (i) holds, and we are done) or $i_0 = 2$ and $r = 2$.
In the latter case $l_1' = m_1' = 0$, $l_2' = l-1$ and $m_2' = m-1$, i.e.,
$(l_1, m_1) = (1, 1)$ and $(l_2, m_2) = (l-1, m-1)$. Condition (b) now
tells us that $l = m = 2$, so that (ii) holds.

This completes the proof of Lemma~\ref{lem.elem} and thus 
of Proposition~\ref{prop.grass1}.
\end{proof}

\section{Proof of Theorem~\ref{thm:GLm} over an algebraically closed field}
\label{sect.GLm}

Recall from Section~\ref{sect.geometric-actions} that $X = \Mnm$ is
the associated $\GL_m \times \PGLn$-variety for the $\GL_m$-action on
$\UD(m, n)$.  Here
$\PGLn$ acts on $\Mnm$ by simultaneous conjugation (since $m \ge 2$,
this action is generically free) and $\GLm$ acts 
on $\Mnm$ by sending $(A_1, \dots, A_m)$
to $(B_1, \dots, B_m)$, where $B_j = \sum_{i = 1}^m
c_{ij} A_i$ and $g^{-1} = (c_{ij})$.  

We shall assume throughout this section that $k$ is an algebraically 
closed field of characteristic zero. Our goal is to prove 
Theorem~\ref{thm:GLm} over such $k$.  In view of \cite[Theorem 1.4(a)]{rv3}, 
we only need to establish the following.

\begin{thm} \label{thm1.grass} 
  Assume $k$ is algebraically closed.
  If $n \geq 3$ and $2 \leq m \leq n^2 - 2$, then the $\GL_m \times
  \PGLn$-action on $\Mnm$ is generically free. 
\end{thm}

\begin{proof}
  The linear action of $\GL_m$ on $\Mnm$ is easily seen to be the
  direct sum of $n^2$ copies of the natural $m$-dimensional
  representation of $\GL_m$, i.e., to be isomorphic to the
  $\GLm$-action on $n^2$-tuples of vectors in $k^m$.  Since $n^2 > m$,
  this action is generically free.  Corollary~\ref{thm.div.cor1}(a)
  with $G = \GL_m$ and $X=\Mnm$ tells us that the $\PGLn$-action on
  $\Mnm/\GLm$ has a reductive stabilizer~$S$ in general position, and that
  $S\simeq\Stab_{\GL_m \times \PGL_n}(x)/ \Stab_{\GL_m}(x)=
  \Stab_{\GL_m \times \PGL_n}(x)$ for $x\in X$ in general position.

  Recall that $\Mnm/\GLm$ is $\PGLn$-equivariantly birationally
  isomorphic to the Grassmannian $\Gr(m, n^2)$ of $m$-dimensional
  subspaces of $\Mn$.  Thus the $\PGL_n$-action on $\Gr(m, n^2)$ has
  a stabilizer $S$ in general position, where $S$ is a reductive
  subgroup of $\PGLn$.  (Recall that $S$ is only well-defined up to
  conjugacy in $\PGLn$).  To prove Theorem~\ref{thm1.grass}, it
  suffices to show that $S$ is trivial.

Assume the contrary. Since $S$ is reductive, 
it contains a non-trivial element $g$ of finite order.  
Then every $L \in \Gr(m, n^2)$ in general position 
is fixed by some conjugate of $g$.  In other
words, the map
\begin{equation}\label{e.dim-count:map}
  \begin{array}{ccc}
    \PGLn \times \Gr(m, n^2)^g  & \lra    & \Gr(m, n^2)\\
                           (h, L)      & \mapsto & h(L)
  \end{array}
\end{equation}
is dominant; here $\Gr(m, n^2)^g$ denotes the fixed points of $g$ in
$\Gr(m,n^2)$.  Denote by $C(g)$ the centralizer of $g$ in $\PGLn$.
Note that $\Gr(m,n^2)^g$ is $C(g)$-stable.  Hence the fiber
of \eqref{e.dim-count:map} over $h(L)$ contains $(hc,c\inv(L))$ for
every $c\in C(g)$.  So by the fiber dimension theorem,
\begin{equation}\label{e.dim-count:new}
    \dim\Gr(m, n^2)+\dim C(g)\leq
    \dim\PGLn+\dim\Gr(m,n^2)^g\,.
\end{equation}
  
  Since $g$ has finite order, it is diagonalizable.  So we may assume that
  \[
     g = \diag( \underbrace{\lambda_1, \dots, \lambda_1}_{\text{$l_1$ times}},
    \dots, \underbrace{\lambda_s, \dots, \lambda_s }_{\text{$l_s$ times}})
    = \diag(\alpha_1,\dots,\alpha_n)
    \, ,
  \]
  where $\lambda_1, \dots, \lambda_s$ are the (distinct) eigenvalues
  of $g$. Note that $s \geq 2$, because $g \neq 1$ in $\PGLn$ and that
  $g$ acts on the matrix units $E_{ij}$ by $g \cdot E_{ij} = \alpha_i
  \alpha_j^{-1} E_{ij}$. The matrix algebra $\Mn$ naturally decomposes
  as a direct sum of character spaces
  \[ V_{\mu} = \Span(E_{ij} \, | \, \alpha_i \alpha_j^{-1} = \mu) \, . \]
  In particular, $V_1$ is the commutator of $g$ in $\Mn$.
  Now~\eqref{e.dim-count:new} implies
  \begin{align*}
    \dim \, \Gr(m, n^2) - \dim \, \Gr(m, n^2)^g
          &\le \dim(\PGLn) - \dim C(g)\\
          &=n^2 - \dim(V_1) \, .
  \end{align*}
  So if $n\geq 3$, part~(b) of the following lemma gives the desired
  contradiction, which completes the proof of
  Theorems~\ref{thm1.grass}, and thus of Theorem~\ref{thm:GLm} in the
  case that $k$ is algebraically closed.
\end{proof}

\begin{lem} \label{lem5.grass}
  Let $n\geq 2$, and $2\leq m\leq n^2 - 2$.
  \begin{itemize}
    \item[\textup{(a)}]$\dim V_1 \ge \dim \, V_{\mu}$ for any $\mu \neq 1$.
    \item[\textup{(b)}]If $n \geq 3$ \textup{(}or $n=2$ and there are
      more than two distinct nonzero $V_\mu$\textup{)}, then $\dim \,
      \Gr(m, n^2) - \dim \, \Gr(m, n^2)^g > n^2 - \dim(V_1)$.
  \end{itemize}
\end{lem}

\begin{proof}
    (a) Note that $\dim \, V_1 = l_1^2 + \dots + l_s^2$ and
  \[ \dim \, V_{\mu} = \sum_{\lambda_i \lambda_j^{-1} = \mu} l_i l_j \]
  Since the eigenvalues $\lambda_1, \dots, \lambda_s$ of $g$ are
  distinct, the last sum has at most one term for each $i = 1, \dots,
  s$. Thus there is a permutation $\sigma$ of $\{ 1, \dots, s \}$ such
  that
  \[ \dim \, V_{\mu} \le  l_1 l_{\sigma(1)} + \dots + l_s l_{\sigma(s)} \, . \]
  So for $v = (l_1, \dots, l_s)$ and $w = (l_{\sigma(1)}, \dots,
  l_{\sigma(s)})$, $\dim V_\mu\leq v\cdot w$.  Hence by
  the Cauchy-Schwarz inequality,
  \[\dim \, V_{\mu} \le v \cdot w \le |v| \, |w| = |v|^2
      = l_1^2 + \dots + l_s^2 = \dim(V_1)\,.\]

    (b) Since $g$ is semisimple, every $L \in \Gr(m, n^2)^g$ is a
    direct sum of its character subspaces spaces, i.e., a direct sum
    of vector subspaces of the $V_{\mu}$. Part (b) now follows from
    Proposition~\ref{prop.grass1} with $V = \M_n$, $N = n^2$, $N_\mu =
    \dim(V_\mu)$, $Z=\Gr(m,n^2)^g$, and $r$ the number of distinct
    nonzero $V_\mu$. 
\end{proof}

\begin{remark} \label{rem.n=2}
We assumed throughout this section that $n \ge 3$. If $n = 2$ then 
the above argument still goes through provided there are more than
two distinct non-zero character subspaces $V_\mu$;
see Lemma~\ref{lem5.grass}(b). In particular, this will always 
be the case if $g^2 \ne 1$ in $\PGL_n$; indeed, in this case
$g = (\lambda_1, \lambda_2)$, where $\mu = \lambda_1/\lambda_2 \ne \pm 1$
and the three spaces $V_1$, $V_{\mu}$ and $V_{\mu^{-1}}$
are distinct. Thus the above argument also shows that for 
$n = m = 2$, either $|S| = 1$ or $S$ has exponent $2$. 
It turns out that, in fact, in this case $|S| = 2$; 
see~\cite[Lemma 14.2]{rv3}.
\end{remark}

\begin{remark} \label{rem.elashvili}
An alternative approach to proving Theorem~\ref{thm1.grass} would be
to appeal to the classification, due to
A. G. Elashvili~\cite{elashvili} and A. M. Popov~\cite{ampopov}, of
pairs $(G, \phi)$, where $G$ is a semisimple algebraic group and $\phi
\colon G \hookrightarrow \GL(V)$ is an irreducible linear
representation of $G$ such that the $G$-action on $V$ has a
non-trivial stabilizer in general position.  Since this classification
is rather involved, and since additional work would be required to apply
it in our situation (note that the group $\GLm \times \PGLn$ is not
semisimple, and that its representation on $\Mnm$ is not irreducible),
we opted instead for the self-contained direct proof presented in this
section.
\end{remark}

\section{$\SLm$-invariant generic matrices}

The goal of this section is to relate the rings of $\SLm$-invariants
in $\Gmn$ and $\UDmn$.

\begin{lem} \label{lem.concom4}\begin{itemize}
  \item[(a)] Every element of $\UDmn^{\GLm}$ can be written in the
    form $\frac{a}{b}$, where $a$ is a homogeneous element of
    $(\Tmn)^{\SLm}$, and $b$ is a non-zero homogeneous element of
    $\Z(\Tmn)^{\SLm}$ of the same degree as~$a$.
  \item[(b)] Assume that a subgroup $G$ of $\GLm$ has no non-trivial
    characters.
    Then every element of $\UDmn^{G}$ can be written as $\frac{a}{b}$
    where $a \in (\Tmn)^{G}$ and $0 \ne b \in \Z(\Tmn)^G$.
\end{itemize}
\end{lem}

\begin{proof} Both parts follows from Lemma~\ref{lem.dc}.
In part (a), we take $\Gamma = \GLm \times \PGLn$, $V = \Mnm$ 
(with the $\Gamma$-action defined in the beginning of Section~\ref{sect.GLm})
and $W = \Mn$ (where $\GLm$ acts trivially on $W$ and $\PGLn$ acts
by conjugation). Here the relative concomitants $\Mnm \lra \Mn$
are the homogeneous elements of $(\Tmn)^{\SLm}$ and 
the relative invariants $\Mnm \lra k$ are the homogeneous elements 
of $\Z(\Tmn)^{\SLm}$; cf.\ Theorem~\ref{thm.procesi}.

If $G$ has no non-trivial characters then relative 
concomitants are (absolute) concomitants, i.e., elements of $(\Tmn)^G$. 
Similarly, relative invariants are elements of $\Z(\Tmn)^G$, and
part (b) is thus simply a restatement of Lemma~\ref{lem.dc} 
in this special case.
\end{proof}
 
\begin{prop} \label{prop.concom5}
Let $G$ be a subgroup of $\GLm$ such that
$G$ has no non-trivial characters.
Then the following conditions are equivalent:
\begin{itemize}
\item[(a)] $\UDmn^G$ has PI-degree $n$.
\item[(b)] $(\Tmn)^G$ has PI-degree $n$.
\item[(c)] $(\Gmn)^G$ has PI-degree $n$.
\end{itemize}
\end{prop}

\begin{proof} The equivalence of (a) and (b) follows from 
Lemma~\ref{lem.concom4}(b). The implication
(c) $\Longrightarrow$ (b) is obvious,
since $\Gmn \subset \Tmn$. It thus remains to prove that 
(b) $\Longrightarrow$ (c).

 Let $g_n$ be the multilinear central polynomial for $n\times
  n$-matrices in \cite[13.5.11]{mcr} (or \cite[p.~26]{rowen:PI}).  If
  $R$ is a prime PI-algebra of PI-degree~$n$, denote by $g_n(R)$ the
  set of all evaluations of $g_n$ in $R$, and denote by $R\,g_n(R)$
  the nonzero ideal of $R$ generated by $g_n(R)$.  Denote by $T$ the
  trace ring of $R$.  (Since we are working in characteristic zero,
  the (noncommutative) trace ring in \cite[13.9.2]{mcr} is the same as
  the one we are using, see \cite[13.9.4]{mcr}.)  Then $R\,g_n(R)$ is
  a common ideal of $R$ and $T$, see \cite[13.9.6]{mcr} (or
  \cite[4.3.1]{rowen:PI}).
 
  Now let $R=\Gmn$.  Then its trace ring is $T = \Tmn$. Recall that 
  we are assuming (b) holds, i.e., $T^G$ has PI-degree~$n$.  Let
  $s$ be a non-zero evaluation of $g_n$ on $T^G$.  Then $s$ is a
  nonzero $G$-invariant, and a central element of $T$ (since it is
  also an evaluation of~$g_n$ on~$T$).  Since $g_n$ is multilinear,
  and since $T$ is generated as an $R$-module by central elements, $s$
  belongs to the ideal of $T$ generated by $g_n(R)$, so that
  $s\,T\subseteq R\,g_n(R)\subseteq R$.  Since $s$ is a
  $G$-invariant, it follows that $s\,T^G \subseteq R^G$.
  Consider the central localization $R^G[s\inv]\subseteq\UDmn$.
  Since it contains $T^G$, $R^G[s\inv]$ must have PI-degree~$n$,
  implying that also $R^G$ must have PI-degree~$n$. This completes
  the proof of the implication (b) $\Longrightarrow$ (c) and thus
  of Proposition~\ref{prop.concom5}.
\end{proof}

\begin{remark} \label{rem:UDmn^GLm:4}
  The same argument also shows that if the three equivalent conditions
  in Proposition~\ref{prop.concom5} are true, then the division
  algebras of fractions of $(\Gmn)^G$ and $(\Tmn)^G$ are both equal to
  $\UDmn^G$.
\end{remark}

\section{Proof of Theorems~\ref{thm.intro1} and~\ref{thm:GLm}}

\subsection*{Proof of Theorem~\ref{thm.intro1}}
Proposition~\ref{prop.concom5} tells us that $(\Gmn)^{\SLm}$ has
PI-degree $n$ if and only if so does $\UDmn^{\SLm}$. Thus
in order to prove Theorem~\ref{thm.intro1} it suffices to show 
that $\UDmn^{\SLm}$ has PI-degree $n$ whenever $2 \le m \le n^2 -2$. 

For $n = m = 2$ we showed this in \cite[Remark 14.4]{rv3} (in fact, 
the argument we gave there remains valid over any base field $k$ 
of characteristic $\ne 2$). For $n \ge 3$ and $2 \le m \le n^2 - 2$,  
Theorem~\ref{thm:GLm} tells us that $\UDmn^{\GLm}$ has PI-degree $n$
(and consequently, so does $\UDmn^{\SLm}$). In summary, we have shown 
that {\em Theorem~\ref{thm.intro1} follows from Theorem~\ref{thm:GLm}}.
\qed

\subsection*{Proof of Theorem~\ref{thm:GLm}}
We have already proved Theorem~\ref{thm:GLm} in the case where
the base field $k$ is algebraically closed; see Section~\ref{sect.GLm}.
We will now reduce the general case to this one by using
Lemma~\ref{lem.concom4}(a).

We begin with a simple lemma. 

\begin{lem} \label{lem.red1}
  Let $K$ be an extension field of $k$, let $V$ be a finite-dimensional
  $k$-vector space, and $V_K = V \otimes_k K$.  Given a linear
  representation of $\SLm(k)$ on $V$, we have
  \[ (V_K)^{\SLm(K)} = V^{\SLm(k)} \otimes_k K \, . \]
\end{lem}

\begin{proof}
  Since $\SLm(k)$ is dense in $\SLm(K)$, the subspace
  $(V_K)^{\SLm(K)}$ is defined inside $V_K$ by a system of homogeneous
  linear equations with coefficients in~$k$.  Clearly finitely many of
  these equations suffice.  Since the dimension of the solution space
  of such a system is the rank of the corresponding matrix (which has
  coefficients in $k$), $(V_K)^{\SLm(K)}$ has a $K$-basis consisting
  of elements of $V^{\SLm(k)}$.
\end{proof}

For the remainder of this section, we will write $\Gmn(K)$, $\Tmn(K)$
and $\UDmn(K)$ to denote the generic matrix algebra, trace ring and
universal division algebra defined over the field $K$.  
Denote the algebraic closure of $k$ by $\overline{k}$. 
Since the
$\SLm$-action on $\UDmn$ preserves degree, Lemma~\ref{lem.red1}
immediately implies the following fact, which we record for later use.

\begin{cor} \label{lem.cor.red1}
  $\Gmn(\overline{k})^{\SLm(\overline{k})}=\Gmn(k)^{\SLm(k)}\otimes_k\overline{k}$,
  and
  $\Tmn(\overline{k})^{\SLm(\overline{k})}=\Tmn(k)^{\SLm(k)}\otimes_k\overline{k}$.
  \qed
\end{cor}

We are now ready to complete the proof of Theorem~\ref{thm:GLm} 
over an arbitrary field $k$ of characteristic zero. 
In Section~\ref{sect.GLm} 
we showed that Theorem~\ref{thm:GLm} holds over the algebraic closure
$\overline{k}$ of $k$. That is, if $2 \le m \le n^2 - 2$ then there
exist elements $c_1, \dots, c_r \in  
\UDmn(\overline{k})^{\GLm(\overline{k})}$ which span
$\UDmn(\overline{k})$ as a vector space over its center.
By Lemma~\ref{lem.concom4} we can write $c_i = a_i/b_i$,
where $a_i \in \Tmn(\overline{k})[d_i]^{\SLm}$ and 
$0 \ne b_i \in \Z(\Tmn(\overline{k}))[d_i]^{\SLm}$ for some $d_i \ge 0$,
$i = 1, \dots, r$.
By Lemma~\ref{lem.red1}, with $K = \overline{k}$ and $V = 
\Z(\Tmn(k))[d_i]$, we have $\Z(\Tmn(k))[d_i]^{\SLm} \ne 0$.
We may now replace $b_i$ by a non-zero element of  
$\Z(\Tmn(k))[d_i]^{\SLm}$. The new $c_i = a_i/b_i$ are still
$\GLm$-invariant elements of $\UDmn(\overline{k})$, and they
still generate $\UDmn(\overline{k})$ as a vector space over its center.

We now apply Lemma~\ref{lem.red1} once again (this time with 
$V = \Tmn(k)[d_i]$) to write each $a_i$ as a finite 
sum $\sum \gamma_{ij} a_{ij}$, where each $\gamma_{ij}\in \overline{k}$ and
each $a_{ij}\in\Tmn(k)[d_i]^\SLm$. 
Now replace our collection of $\GLm$-invariant
elements $\{ c_i = a_i/b_i \}$ in $\UDmn(\overline{k})$
by $\{ c_{ij} = a_{ij}/b_i \}$.
By construction, the elements $c_{ij}$ lie in $\UDmn(k)^{\GLm}$
and span $\UDmn(\overline{k})$ over its center. Hence,
these elements generate a $k$-subalgebra of $\UDmn(k)^{\GLm}$
of PI-degree $n$. Consequently, $\UDmn(k)^{\GLm}$ itself 
has PI-degree $n$. This completes the proof of 
Theorem~\ref{thm:GLm} (and of Theorem~\ref{thm.intro1}).
\qed

\section{The case $m \ge n^2 - 1$}
\label{sect.m>=n^2-1}

Theorems~\ref{thm.intro1} and~\ref{thm:GLm} assume that $m \le n^2 -
2$.  We will now describe $\UDmn^{\GLm}$ and $\UDmn^{\SLm}$ for $m \ge
n^2 -1$.

Recall the definition of the discriminant of $n^2$ matrices of size
$n\times n$, say $A_1,\ldots,A_{n^2}$: it is the determinant of the
$n^2\times n^2$-matrix whose $i$-th row consists of the entries of
$A_i$, cf.\ \eqref{eq.n^2 x n^2}.  When viewed as a function
$(\Mn)^{n^2} \lra k$, $\Delta$ is the unique multilinear alternating
function such that $\Delta(e_{11}, e_{12}, \dots, e_{nn}) = 1$; cf.,
e.g., \cite[Lemma 3]{formanek4}.  Here the $e_{ij}$ are the matrix
units.

\begin{prop} \label{prop.m=n^2}
{\upshape(a)} If $m > n^2$, then $\UDmn^{\SLm} = \UDmn^{\GLm} = k$. 

\smallskip

\noindent
Now let $m = n^2$, and denote by $\Delta$ the discriminant of the
generic matrices $X_1,
\dots, X_m$.
\begin{itemize}
\item[(b)] $\UDmn^{\GLm} = k$.
\item[(c)] $(\Tmn)^{\SLm} = k[\Delta]$.
\item[(d)] $\UDmn^{\SLm} = k(\Delta)$.
\end{itemize}
\end{prop}

\begin{proof}
(a) We may clearly assume that $k$ is algebraically closed.  In this
case $\SLm$ has a dense orbit in  
the associated variety $X = \Mnm$.  Thus the rational quotient $X/\SLm$
is a single point (with trivial $\PGLn$-action), and
\[ \UDmn^{\SLm} = \RMaps_{\PGLn}(pt, \Mn) = k \, . \]
Now suppose $m = n^2$. Then $\GLm$ has a dense orbit in $X = \Mnm$. Arguing 
as in part (a), we prove (b); cf.~\cite[Proposition 13.1(a)]{rv3}.
(c) is proved in~\cite[p. 210]{formanek3}, and (d) follows from (c)
by Lemma~\ref{lem.concom4}(b).
\end{proof}

\begin{remark}\label{rem.formanek.m=n^2}
  Let $m=n^2$.  Formanek showed that $\Delta \not \in (\Gmn)^\SLm$
  (\cite[p.~214]{formanek3}) but $\Delta^i \in \Gmn$ for every integer
  $i \ge 2$ (this follows from \cite[Theorem~16]{formanek4}).
  Consequently for $m=n^2$,
  \[(\Gmn)^\SLm=k[\Delta^2,\Delta^3]\,.\]
\end{remark}

\begin{prop} \label{prop.m=n^2-1} 
Suppose $m = n^2-1$, and let 
\[Y = \sum_{i, j = 1}^n\Delta(X_1, \dots, X_m, e_{ji})\, e_{ij}\,,\]
where the $e_{ij}$ are the matrix units.
\begin{itemize}
\item[(a)] $Y \in (\Tmn)^{\SLm}$.
\item[(b)] The eigenvalues of $Y$ are algebraically independent over
  $k$ {\upshape(}and, in particular, distinct{\upshape)}.
\item[(c)] $(\Tmn)^{\SLm} = k[c_1, \dots, c_{n-1}, Y]$ is a
  polynomial ring in $n$ independent variables over $k$.  Here $c_1
  = - \tr(Y), \dots, c_n = (-1)^n \det(Y)$.
\item[(d)] $\UDmn^{\SLm} = k\bigl(c_1, \dots, c_{n-1}, Y\bigr)$.
\item[(e)] $\UDmn^{\GLm} = k\Bigl(\frac{c_2}{(c_1)^2}, \dots,
    \frac{c_{n-1}}{(c_1)^{n-1}}, \frac{1}{c_1}Y\Bigr)$.
\end{itemize}
\end{prop}

\begin{proof} 
For the proof of (a)---(c), we may assume that $k$ is algebraically
closed (cf.\ Corollary~\ref{lem.cor.red1}).
(a) We view $Y$ as a regular map $\Mnm \lra \Mn$.
We want to show that this map is $\PGLn$-equivariant, i.e., $Y \in \Tmn$.
Since $Y$ is clearly $\SLm$-equivariant (recall that $\SLm$
acts trivially on $\Mn$), this will imply part~(a).

We begin by observing that for any
$(A_1, \dots, A_m) \in \Mnm$, and $Z \in \Mn$, 
\begin{equation} \label{e.Y}
\tr(Y(A_1, \dots, A_m) Z) = \Delta(A_1, \dots, A_m, Z) \, . 
\end{equation}
Indeed, both sides are linear in $Z$, so we only need to
check~\eqref{e.Y} for the elementary matrices $Z = e_{ij}$, 
where it is easy to do directly from the definition of $Y$.

\smallskip Fix an $m$-tuple $(A_1, \dots, A_m) \in \Mnm$ of $n \times
n$-matrices.  Since the trace form on $\Mn$ is non-singular, $Y(A_1,
\dots, A_m)$ is the unique matrix satisfying~\eqref{e.Y} for every $Z
\in \Mn$. The $\PGLn$-equivariance of $Y \colon \Mnm \lra \Mn$ is an
easy consequence of this and the fact that $\Delta$ is
$\PGLn$-invariant (see \cite[p.~209]{formanek3}). This concludes the
proof of (a).

\smallskip
Our proof of parts (b) and (c) relies on the following claim:
$Y \colon \Mnm \lra \Mn$ is the categorical quotient 
map for the $\SLm$-action on $\Mnm$. In other words, we claim 
that the $n^2$ elements $\Delta(X_1, \dots, X_m, e_{ij})$ 
($i, j = 1, \dots, n$) generate $k[\Mnm]^{\SLm}$ as a $k$-algebra.
To prove this claim we will temporarily write 
$(A_1, \dots, A_m) \in \Mnm$ as an $m \times n^2$-matrix
\begin{equation}\label{eq.n^2 x n^2} 
  A = \begin{pmatrix} a_{11}^{(1)} & a_{12}^{(1)} & \dots &a_{ij}^{(1)} & \dots & a_{nn}^{(1)} \\ 
                    \vdots         & \vdots       &       & \vdots      &       & \vdots       \\
                   a_{11}^{(m)} & a_{12}^{(m)} & \dots &a_{ij}^{(m)} & \dots & a_{nn}^{(m)}
  \end{pmatrix} \, . 
\end{equation}
That is, we write each $n \times n$ matrix $A_h = (a_{ij}^{(h)})$ as a
single row of $A$. In this notation, $g\in\SLm$ acts on $\Mnm$ by
multiplication by the transpose of $g\inv$ on the left; that is, $g(A)
= (g\inv)^{\text{transpose}} \cdot A$ for every $g \in \SLm$. Let
$\delta_{ij}(A_1,\ldots,A_m)$ be the $m \times m$-minor of this matrix
obtained by removing the $ij$-column from $A$ and taking the
determinant of the resulting $m \times m$-matrix. The first 
theorem of classical
invariant theory (see~\cite[Theorem II.6.A]{weyl} or
\cite[Theorem~2.1]{dolgachev}) says that the elements
$\delta_{ij}(X_1,\ldots,X_m)$ generate $k[\Mnm]^{\SLm}$ as a
$k$-algebra.  On the other hand, it is easy to see that
$\delta_{ij}(X_1,\ldots,X_m) = \pm \Delta(X_1, \dots, X_m, e_{ij})$.
This proves the claim.

Now observe that since $m=n^2-1$,
\[ \dim(\Mnm \catq \SLm) = mn^2 - (m^2-1) = n^2 = \dim(\Mn) \, . \]
This means that the $n^2$ $\SLm$-invariant functions
\[ \Delta(X_1, \dots, X_m, e_{ji}) \colon \Mnm \lra k \]
are algebraically independent over $k$. In other words, $Y$ (viewed as
a matrix in $\Tmn \subset \Mn(k[x_{ij}^{(h)}]$) has algebraically 
independent entries. Part (b) easily follows from this assertion;
cf.~\cite[Lemma II.1.4]{procesi1}. 

\smallskip
Furthermore,
\begin{align*} (\Tmn)^{\SLm} &= \Morph_{\SLm \times \PGLn}(\Mnm,\Mn)\\
    & \simeq \Morph_{\PGLn}(\Mnm \catq \SLm, \Mn)\\
    & \simeq\Morph_{\PGLn}(\Mn, \Mn) = T_{1, n} \, , 
\end{align*}
where $T_{1, n}$ is the trace ring of one generic $n \times n$-matrix.
Here the last equality is a special case of 
Procesi's Theorem~\ref{thm.procesi}(a) (with $m = 1$).
Since the chain of isomorphisms identifies $Y$ with the identity map
$\Mn \lra \Mn$, we conclude that
\[ (\Tmn)^{\SLm} = k[c_1, \dots, c_n, Y]  \, . \]
Since $ Y^n + c_1 Y^{n-1} + \dots + c_n = 0$,
$ k[c_1, \dots, c_{n-1}, Y] = k[c_1, \dots, c_n, Y]$.
This proves the first assertion in (c).

To show that $c_1, \dots, c_{n-1}, Y$ are algebraically 
independent over $k$, denote the eigenvalues of $Y$ by 
$\lambda_1, \dots, \lambda_n$.  By part (b), $\lambda_1, \dots, 
\lambda_n$ are algebraically independent over $k$.  Since 
$Y$ is algebraic over $k(c_1, \dots, c_n)$, we have 
\begin{align*} \trdeg_k &\,k(c_1, \dots, c_{n-1}, Y) = \trdeg_k
k(c_1, \dots, c_n, Y) \\
&=\trdeg_k k(c_1, \dots, c_n) =
\trdeg_k k(\lambda_1, \dots, \lambda_n) = n \, . 
\end{align*}
This shows that $c_1, \dots, c_{n-1}, Y$ are algebraically 
independent over $k$, thus completing the proof of (c).

(d) is an immediate consequence of (c) and Lemma~\ref{lem.concom4}(b).
To prove (e), denote the central torus of $\GLm$ by $\bbG_m$. Then 
\[ \UDmn^{\GLm} = (\UDmn^{\SLm})^{\bbG_m} = k(c_1, \dots, c_{n-1}, Y)^{\bbG_m}
\, , \]
where $\bbG_m$ acts on the purely transcendental extension 
$k(c_1, \dots, c_{n-1}, Y)$ as follows: 
$t \cdot c_i \mapsto t^{im} c_i$ for $i =1, \dots, n-1$, and
and $t \cdot Y \mapsto t^m Y$. Part (e) easily follows from 
this description.
\end{proof}

\begin{remark}\label{rem.formanek.m=n^2-1}
Note that $c_1 = - \Delta(X_1,\ldots,X_{n^2-1}, I_n)$,
where $I_n$ is the $n\times n$ identity matrix.   
By a theorem of Formanek, $(c_1)^2$
is an element of $\Z(\Gmn)^{\SLm}$ for $m = n^2 -1$, see
\cite[Theorem~16]{formanek4}.
\end{remark}

\section{Proof of Theorem~\ref{thm.intro2}}

By Corollary~\ref{lem.cor.red1}, we may assume that $k$ is
algebraically closed.  Set $A=(\Gmn)^\SLm$ and $B=(\Tmn)^\SLm$.  By
Theorem~\ref{thm.intro1}, $A$ and $B$ both have PI-degree~$n$.  Thus
$\Z(A)=(\Z(\Gmn))^\SLm$ and $\Z(B)=(\Z(\Tmn))^\SLm$.
Since $\SLm$ is a reductive group, and since
$\Tmn$ is a finitely generated $k$-algebra and a finite module over
its center, $B$ is a finite $\Z(B)$-module, and both $B$ and $\Z(B)$
are finitely generated Noetherian $k$-algebras, cf.\ 
\cite[Proposition~4.2]{vonessen:memoir}.  Moreover, $B$ is an
FBN ring, cf.\ \cite[13.6.6]{mcr}.

By Corollary~\ref{thm.div.cor1}(b) and Remark~\ref{rem:UDmn^GLm:4},
the transcendence degrees of both $B$ and $\Z(B)$ are
$t=(m-1)n^2-m^2+2$. For notational simplicity, set  
 \[\mu(S)= \limsup_{d\to\infty}\,\frac{\dim_k S[d]}{d^{t-1}} \]
for any graded $k$-algebra $S = \oplus_{d \ge 0} S[d]$.
By \cite[Lemma 6.1]{stafford-zhang} (cf.\ also
\cite[\S12.6]{kl}), $f(d)=\dim_k B[d]$ is eventually periodically
polynomial, i.e., there are polynomials $f_1,\ldots, f_s$ with
rational coefficients such that $f(d)=f_i(d)$ whenever $d$ is large
enough and congruent to $i$ modulo $s$; moreover, the maximum of the
degrees of the $f_i$ is $t-1$.  Consequently $\mu(B)$
exists and is equal to the largest among the leading coefficients of
those $f_i$ of degree $t-1$.  A similar argument shows that
$\mu(Z(B))$ exists and is a nonzero number.

Consider the multilinear central polynomial $g_n$ for $n\times n$
matrices used in the proof of Proposition~\ref{prop.concom5}.  Since
it is multilinear and nonzero on $A$, we can find a nonzero
evaluation~$c$ of $g_n$ at homogeneous elements of $A$; this $c$ is
homogeneous.  Since $c$ is also an evaluation of $g_n$ on $\Gmn$,
$c\Tmn\subset \Gmn$, so that $cB\subset A$ and
$c\Z(B)\subset \Z(A)$.  Then for all integers $d \ge j$,
$cB[d-j]\subseteq A[d]\subseteq B[d]$, where $j=\deg c$.  Replacing
$c$ by $c^s$ if necessary, we may assume that $s$ divides $j$.
Consequently, whenever $d$ is large enough and congruent to $i$ modulo
$s$,
\[f_i(d-j)\leq \dim_k A[d]\leq f_i(d)\,.\]
It follows easily that $\mu(A)$ exists and is equal 
to the largest among the leading coefficients of
those $f_i$ of degree $t-1$.  A similar argument shows that
$\mu(Z(A))$ exists and is a nonzero number.
\qed 

\begin{remark} 
  The above proof shows that
  $\mu((\Gmn)^\SLm)=\mu((\Tmn)^\SLm)$, and that
  $\mu(\Z((\Gmn)^\SLm))= \mu(\Z((\Tmn)^\SLm))$.
\end{remark}

\begin{remark} \label{rem9.2} 
Consider the $\GLm$-representation on $R$, where
$R = \Gmn$, $\Tmn$, $\Z(\Gmn)$ or $\Z(\Tmn)$.
Recall that irreducible polynomial representations 
of $\GLm$ are indexed by partitions $\lambda = (\lambda_1,
\dots, \lambda_s)$ with $s \le m$ parts; 
cf.~\cite[Section 2]{formanek3}. 
Denote the multiplicity of the irreducible representation 
corresponding to $\lambda$ in $R$ by $\mult_{\lambda}(R)$.
Writing $(r^m)$ for the partition 
$\lambda = (r, ..., r)$ ($m$ times), we have

\smallskip
(a) $ \dim_k \, R^{\SLm}[d] = 0$ if $d$ is not a multiple of $m$, and

\smallskip
(b) $ \dim_k \, R^{\SLm}[rm] = \mult_{(r^m)}(R)$ for any integer $r \ge 1$.
\end{remark} 

\begin{proof}
  (a) Assume $R^\SLm[d]$ is nonzero.  Then it is a direct sum of
  one-dimensional representations of $\GLm$ of the form $M =
  \Span_k(v)$.  Moreover, any such representation is given by $g(v) =
  \det(g)^r v$ for some integer $r$; cf., e.g., \cite[Theorem
  3(a)]{formanek3}.  On the other hand, substituting $g = tI_m$, where
  $t \in k$ and $I_m$ is the $m \times m$ identity matrix, we obtain,
  $g(v) = t^d v$.  Since $\det(t I_m) = t^m$, we see that $d = rm$, as
  claimed.

\smallskip
(b) If $d = rm$ and $0 \ne v \in R^{\SLm}[rm]$
then the partition associated to the 1-dimensional 
irreducible $\GLm$-module $M = \Span(v)$ 
is $(r^m)$; cf. , e.g., \cite[Theorem 2]{formanek3}.
Now consider the direct sum decomposition $R =\oplus R_\lambda$,
where $R_\lambda$ is the sum of all irreducible 
$\GLm$-submodules of $R$ with associated 
partition $\lambda$. The argument of 
part (a) shows that $R_{(r^m)} = R^{\SLm}[rm]$.
Moreover, since $\dim(M) = 1$, we have
\[ \dim_k \, R^{\SLm}[rm] = \dim_k \,  R_{(r^m)} = \mult_{(r^m)}(R) \, , \]
as claimed.
\end{proof}

\section{Standard polynomials}
Let $\Gmn$ be the ring of $m$ generic $n \times n$-matrices.
By Theorem~\ref{thm.intro1}, $(\Gmn)^{\SL_m}$ is a PI domain 
of degree $n$, whenever $2 \le m \le n^2 -2$. We will now
describe one particular element of this ring. Let
\[ F_m(x_1, \dots, x_m) = 
\sum_{\sigma \in \Sym_n} (-1)^{\sigma} x_{\sigma(1)} \dots x_{\sigma(n)} 
 \in k \{ x_1, \dots, x_m \} \]  
be the standard polynomial.
Since $F_m$ is multilinear and alternating,
one checks easily that for $g\in\GLm$, 
\begin{equation} \label{e.standard}
g(F_m)=\det(g)\cdot F_m;
\end{equation}
see, e.g., \cite[1.4.12]{rowen:PI}.  
Substituting $m$ generic $n \times n$-matrices $X_1, \dots, X_m$
into $F_m$, we obtain  
\[ f_{m, n} = F_m(X_1, \dots, X_m) = 
\sum_{\sigma \in \Sym_n} (-1)^{\sigma} X_{\sigma(1)} \dots X_{\sigma(n)} 
\in \Gmn \, . \]  
From~\eqref{e.standard}, we see that $f_{m, n} \in (\Gmn)^{\SL_m}$.
By the Amitsur-Levitzki Theorem, $f_{m, n} = 0$ iff $m\geq 2n$.

Fix $m, n \ge 2$ and let $K$ be the center of $\UDmn$. 

\begin{prop}\label{lem:Fm}
  For $2\leq m<2n$, $K(f_{m, n})$ generates a $\GLm$-stable
  maximal subfield of $\UDmn$.
\end{prop}

The proof is algebraic in nature and works in characteristic $\neq 2$.

\begin{proof} The fact that $K(f_{m, n})$ is a $\GL_m$-stable 
subfield follows from~\eqref{e.standard}.
 In order to prove that this subfield is maximal, it suffices to
  verify that $f_{m, n}$ has an eigenvalue of multiplicity $1$.
(Indeed, if, say, $n=d\cdot[K(f_{m, n}):K]$, then the characteristic 
polynomial $p(t)$ of $f_{m, n}$ in $\UDmn$ has the form $p(t) = q(t)^d$, 
where $q(t)$ is the minimal polynomial of $f_{m, n}$ over
$K$. This shows that the multiplicity of each eigenvalue 
of $f_{m, n}$ is divisible by $d$.)

  Since the multiplicity of eigenvalues cannot decrease 
  when evaluating $f_{m, n}$ in $\Mn$,
  it suffices now to show that $f_{m, n}$ (or equivalently, $F_m$)
  has some evaluation in $\Mn$ with an eigenvalue of multiplicity one.
  We now proceed to construct such an evaluation.  Since 
\[ F_m(1,x_2,\ldots,x_m)=F_{m-1}(x_2,\ldots,x_m) \]
 for $m$ odd (cf.~\cite[Exercise 1.2.3]{rowen:PI}), 
 we may assume that $m$ is even, say $m=2r-2$,
  with $1<r\leq n$.  In $\Mn$, consider the sequence of $m$ matrix
  units
  \[e_{1,2},e_{2,2},e_{2,3},e_{3,3},\ldots,
      e_{r-2,r-1},e_{r-1,r-1},e_{r-1,r},e_{r,1}\,.\]
  When permuting these matrix units cyclically, their product is
  nonzero; for any other permutation, their product is zero.  Since
  an $m$-cycle is odd, it follows that $F_m$ evaluated at these matrix
  units is
  \[e_{1,1}-e_{2,2}+e_{2,2}-+
      \cdots-e_{r-1,r-1}+e_{r-1,r-1}-e_{r,r} = e_{1,1}-e_{r,r}\,,\]
  which has $1$ as an eigenvalue of multiplicity one (since
  $\Char(k)\neq2$).
\end{proof}

We do not know an explicit expression for any non-constant element of
$(\Gmn)^{\SL_m}$ (as a polynomial in the generic $n \times n$-matrices
$X_1, \dots, X_m$) in the case where $2n \le m \le n^2 -2$; we leave
this as an open question. 
Note that for $m=n^2$ and $m=n^2-1$,
such elements are exhibited in Remarks~\ref{rem.formanek.m=n^2}
and~\ref{rem.formanek.m=n^2-1}.

\section*{Acknowledgments}
We are grateful to A.\ Berele and E.\ Formanek for helpful comments.

\end{document}